\documentclass{jhrs}

\usepackage{amsmath,amssymb,amscd,amsthm,mathdots,indentfirst}
\usepackage{amsfonts}
\usepackage[mathscr]{eucal}
\usepackage{mathpazo}
\usepackage{ifpdf,ifdraft}

\ifpdf
\usepackage[pdftex,bookmarks,colorlinks,breaklinks]{hyperref}  
\hypersetup{linkcolor=red,citecolor=blue,filecolor=magenta,urlcolor=blue} 
\hypersetup{linkcolor=black,citecolor=black,filecolor=black,urlcolor=black} 
\fi

\newtheorem{prop}[equation]{Proposition}
\newtheorem{thm}[equation]{Theorem}
\newtheorem{cor}[equation]{Corollary}
\newtheorem{lem}[equation]{Lemma}

\theoremstyle{definition}
\newtheorem{defn}[equation]{Definition}
\newtheorem{nota}[equation]{Notation}
\newtheorem{example}[equation]{Example}
\newtheorem{remark}[equation]{Remark}

\numberwithin{equation}{section}
\renewcommand{\theequation}{\thesection.\arabic{equation}}

\def\map{\operatorname{map}}
\def\Map{\operatorname{Map}}

\def\vp{\varphi}
\def\al{\alpha}
\def\be{\beta}

\def\la{\lambda}

\def\De{\Delta}
\def\si{\sigma}
\def\th{\theta}

\def\B{\ensuremath{\mathcal{B}}}
\def\C{\ensuremath{\mathcal{C}}}
\def\D{\ensuremath{\mathcal{D}}}
\def\E{\ensuremath{\mathcal{E}}}

\def\O{\ensuremath{\mathcal{O}}}

\def\fff{\ensuremath{\mathbf{f}}}
\def\ggg{\ensuremath{\mathbf{g}}}
\def\hhh{\ensuremath{\mathbf{h}}}

\def\ev{\operatorname{ev}}

\def\Top{\mathbf{Spc}}

\def\END{\operatorname{END}}

\def\Nat{\operatorname{Nat}}
\def\opEND{\operatorname{\O End}}

\def\Hom{\mathrm{Hom}}

\DeclareMathOperator{\id}{id}
\def\Id{\textrm{Id}}
\def\op{\ensuremath{\mathrm{op}}}
\def\hocolim{\ensuremath{\operatornamewithlimits{hocolim}}}
\def\holim{\ensuremath{\operatornamewithlimits{holim}}}

\def\Lan{\operatornamewithlimits{LKan}}

\def\Sing{\operatornamewithlimits{Sing}}
\def\colim{\operatorname{colim}}

\def\sets{\mathbf{Sets}}

\def\assoc{\operatorname{assoc}}

\def\sk{\operatorname{sk}}
\def\Tot{\operatorname{Tot}}

\def\pnt{{\{ * \} }}

\def\colim{\operatorname{colim}}

\newcommand{\xto}[1]{\xrightarrow{#1}}

\newcommand{\hhocolim}[1]{\operatornamewithlimits{hocolim}_{#1}}
\newcommand{\hholim}[1]{\operatornamewithlimits{holim}_{#1}}

\DeclareMathOperator{\im}{im}

\newcommand{\comp}{\mathbin{\raise0.2ex\hbox{\tiny$\diamond$}}}
\newcommand{\bigcomp}{\mathop{\Diamond}}
\newcommand{\extcomp}{\mathbin{\tilde {\raise0.2ex\hbox{\tiny$\diamond$}}}}
\newcommand{\bigextcomp}{\mathop{\widetilde{\vphantom{J}\Diamond}}}
\newcommand{\exttimes}{\tilde \times}
\newcommand{\bigexttimes}{\mathop{\widetilde\prod}}

\ifdraft{

}

\ifpdf
\else
\def\texorpdfstring#1#2{#1}
\fi

\input xy
\xyoption{all}

\def\smashop#1_#2{%
\displaystyle{#1_{%
\hbox to 0pt{\hss$\scriptstyle{#2}$\hss}}\;}}

\def\smashopsup#1^#2{%
\displaystyle{#1^{%
\hbox to 0pt{\hss$\scriptstyle{#2}$\hss}}\;}}

\def\Deplus{\Delta_+}
\def\De{\Delta}
\def\Dplus{D_+}

\def\D{D}
\def\TED{\widetilde{ED}} 

\def\AHP{\tilde A}

\date\today
\title{$A_\infty$-monads and completion}
\author{Tilman Bauer}
\email{tilman@few.vu.nl}
\address{Faculteit der exacte wetenschappen\\ Vrije Universiteit Amsterdam\\De Boelelaan 1081A\\1081HV Amsterdam, The Netherlands}

\author{Assaf Libman}
\email{assaf@maths.abdn.ac.uk}
\address{
	Department of Mathematical Sciences, 
	King's College\\
	University of Aberdeen\\
  	Aberdeen AB24 3UE , Scotland, U.K.}


\begin{document}

\begin{abstract}
Given an operad $A$ of topological spaces, we consider $A$-monads in a topological category $\C$.
When $A$ is an $A_\infty$-operad, any $A$-monad $K \colon \C \to \C$ can be thought of as a monad up to coherent homotopies.
We define the completion functor with respect to an $A_\infty$-monad and prove that it is an $A_\infty$-monad itself.
\end{abstract}

\maketitle

\section{Introduction}
The starting point for this paper is the claim made in \cite[I.5.6]{yellowmonster} that the $R$-completion functor
$\hat{R}$, where $R$ is a commutative ring, is a monad (triple).
This claim was not proved and later retracted by the authors. 
It is a monad up to homotopy, as was proved in \cite{libman:homotopy-limits-triples,bousfield:cosimplicial-resolutions}. 
In this paper we prove that these homotopies are part of a system of higher coherent homotopies, turning $\hat{R}$ into what we call an \emph{$A_\infty$-monad}. In fact we work in a wider context.

Let $\E$ be a category which is complete and cocomplete, tensored and cotensored over a ``convenient'' category $\Top$ of topological spaces in the sense of Steenrod \cite{steenrod:67} -- see Section \ref{s:enriched-categories}. 
Assume further that $\E$ has an enriched monoidal structure $\comp$ with unit $I$.

For a (non-symmetric) operad $A$ of topological spaces, an \emph{$A$-algebra} in $(\E,\comp,I)$ is an object $K \in \E$ which is equipped with appropriate $A$-algebra structure maps $A(n) \to \map_{\E}(K\comp\dots\comp K,K)$. If $A$ is the associative operad $\assoc$ with $\assoc(n) = *$ for all $n$, an $A$-algebra is just a monoid in $(\E,\comp,I)$ and therefore has an associated cobar construction $R_K\colon \Delta \to \E$, $R_K(n) = K^{\comp (n+1)}$.
Its totalization $\hat K = \Tot R_K$, where $\Delta$ is the category of ordered sets $[n]=\{0,\dots,n\}$ and monotonic maps, is called the \emph{completion} of $K$. 
This construction can be generalized to arbitrary $A$-algebras, as we will show in Section~\ref{sec-a-monads-completion} below. 
The main result of this paper is the following completion theorem:

\begin{thm}\label{completion-theorem}
For any operad $A$ there exists an operad $\hat{A}$ mapping to $A$ with the following property: if $K$ is an $A$-algebra in $(\E,\comp,I)$ then $\hat{K}$ is an $\hat{A}$-algebra and $\hat{K} \to K$ is a morphism of $\hat{A}$-algebras. If $A$ is an $A_\infty$-operad then so is $\hat{A}$.
\end{thm}

Our principal application of this theorem is when $\E$ is the (large) category of continuous endofunctors of a complete and cocomplete category $\C$ tensored and cotensored over $\Top$, with the monoidal structure given by composition of functors.  See Section \ref{ss:endofunctorCat} below for more details, particularly about the set-theoretic difficulties involved.
A monoid $K$ in this category is simply a monad, and an $A$-algebra is what we call an \emph{$A$-monad} on $\C$. 
If the spaces of $A$ are weakly contractible we  call $K$ an $A_\infty$-monad.
The resolution $R_K$ and the completion $\hat K$ generalize the classical definitions given by Bousfield and Kan in \cite{yellowmonster} for monads associated to commutative rings. 
Thus we obtain as an immediate corollary of the theorem:

\begin{cor}\label{completion-corollary}
Let $K$ be an $A$-monad in $\C$.
Then $\hat{K}$ is an $\hat{A}$-monad and $\hat K \to K$ is a morphism of $\hat{A}$-monads.
If $K$ is an $A_\infty$-monad, then so is $\hat K$.
\end{cor}

For example, the theorem works in the following setting:
\begin{itemize}
\item The functor $K(X) = \Omega^\infty(E \wedge X)$ for an $A_\infty$-ring spectrum $E$. If $E$ is ordinary homology with coefficients in a ring $R$ then $\hat K$ is a topological version of the Bousfield-Kan $R$-completion. For more general $E$, it is the Bendersky-Thompson completion \cite{bendersky-thompson}.
\item 
The functor $\hat K$ for an $A_\infty$-monad $K$. The output of the theorem, namely $\hat K$, can be fed back into it to obtain $\hat{\hat{K}}$.
In particular, if we start with a completion at a homology theory $E$ as above, we get a secondary multiplicative Bousfield spectral sequence whose $E_1$-term consists of the homotopy groups of iterated completions of a space, and which converges to the homotopy groups of the completion at $\hat K$, namely $\hat{\hat{K}}(-)$

Following Dror and Dwyer in \cite{dror-dwyer:long-localization}, the completion procedure can be iterated transfinitely to obtain a ``long tower'' $\dots \to \hat{\hat{K}} \to \hat{K} \to K$.
When $E$ is the ordinary homology with coefficients in $\mathbf F_p$ or a subring of $\mathbf Q$, they have shown that if we iterate this often enough, the spectral sequence will eventually converge to the homotopy groups of $L_p X$, the localization. It would be very interesting to know if this works for any homology theory $E$.

\item 
The same construction works in the category of spectra, where the monad is given by $K(X) = E \wedge X$ for an $A_\infty$-spectrum $E$. The theorem allows us to construct a secondary Adams spectral sequence with $E^1$-term constructed from the homotopy groups of the iterated $E$-nilpotent completion and converging to the homotopy groups of a second-order completion. In this way, just as for spaces, a long tower of spectra under the $E$-localization can be constructed, with a higher Adams spectral sequence relating every tower stage to the following one.
\end{itemize}

The crux of the theorem is the assertion that $\hat{A}$ is an $A_\infty$-operad if $A$ is one.
Thus, the completion of a monad is not a monad, but the completion of a monad up to coherent homotopies turns out to be a monad of the same type.
We remark that if $\E$ is enriched over $s\sets$, the category of simplicial sets, the construction of the 
operad $\hat{A}$ and its action on $\hat{K}$ in Sections \ref{sec-a-monads-completion}--\ref{restriced-operad-section} work without any change.
However, it is not true that $\hat{A}$ is an $A_\infty$-operad of simplicial sets if $A$ has the same property (see 
Theorem~\ref{thm-hatA-contractible}.)
For this reason we will only be interested in enrichment over $\Top$.

In the companion paper \cite{bauer-libman:monads2} we show that a related result, stronger in some sense and weaker in another, holds in the simplicial world. We show there that if the monad under consideration is the free $R$-module functor for a ring $R$ then there is an explicit, combinatorially defined $A_\infty$-operad action on the cobar construction, not just its totalization.

It is worth noting that when $(\E,\comp,I)$ is $(\Top,\times,*)$ or any other monoidal category where $I$ is a terminal object then for any $A_\infty$-algebra $K$, $R_K$ has a cosimplicial retraction to $I$ and thus $\hat K = I$. The theorem is thus only making a trivial statement.

\section{Enriched homotopy limits}
\label{s:enriched-categories}

The purpose of this section is to set up some notation and recall various easy facts about enriched categories. Most of the facts can be found in the literature, e.g. \cite{borceux:voltwo,kelly:enriched}, and are collected here for the reader's convenience.

Let $\Top$ denote a complete and cocomplete ``convenient'' category of topological spaces in the sense of Steenrod \cite{steenrod:67}. 
Thus, $\Top$ has mapping spaces with the right adjunction properties where $\map(X,Y)$ is the set of all the continuous functions $X \to Y$ topologized appropriately.
In addition, $\Top$ is large enough to contain the subcategory of CW-complexes and it is equipped with a model structure where fibrations are Serre fibrations and weak equivalences are weak homotopy equivalences. 
The prime examples for $\Top$ are the category of compactly generated Hausdorff spaces \cite{steenrod:67} or the category of $k$-spaces, \cite[Appendix]{lewis:thesis} or \cite[Section 2.4]{hovey:model-categories}. 

Let $C$ be a small category enriched over $\Top$ (in short: a $\Top$-category). Following \cite{hollender-vogt}, we make the following definition.

\begin{defn}\label{toplogicalcategories}
A \emph{CW-category} $C$ is a $\Top$-category all of whose morphism spaces $\map_C(x,y)$ are CW-complexes, and where the identity morphisms are $0$-cells. The composition maps are required to be cellular maps.\end{defn}

Any $\Top$-category $C$ can be
turned into a CW-category by applying the singular and the realization functor
(cf. \cite[VII.\S 2]{yellowmonster}) to its morphism spaces.

We write $C^*$ for the $\Top$-category of continuous (i.e. $\Top$-enriched) functors $C \to \Top$, 
cf. \cite{kelly:enriched}.
This category has all limits and colimits because they are formed in $\Top$. Consequently, the \emph{coend} $F \otimes_C G$ of $F \in {C^{\op}}^*$ and $G \in C^*$ and the \emph{end} 
$[F,G]_C$ of $F,G \in C^*$ are defined \cite[\S 2.1 eq. (2.2) and \S 3.10 eqs. (3.68, 3.70)]{kelly:enriched}.

Every functor $J \colon C \to D$ of small CW-categories induces a continuous restriction functor 
$J^* \colon D^* \to C^*$.
The \emph{left Kan extension} of $F \in C^*$ along $J$ is given by
\[
\Lan_JF = \map_D(J(-), -) \otimes_C F = \int\limits^{c \in C} \map_D(J(c),-) \times F(c) \quad\text{\cite[(4.25)]{kelly:enriched}.}
\]
It is left adjoint to $J^*$, i.e. for every $G \in D^*$ there are natural homeomorphisms
\begin{equation}\label{lkanadjoint}
\map_{D^*}(\Lan_J F,G) \cong \map_{C^*}(F,J^*G) \quad \text{\cite[Theorem~4.38]{kelly:enriched}.}
\end{equation}
It follows that the assignment $F \mapsto \Lan_J F$ gives rise to a continuous functor
\begin{equation}\label{lancontinuous}
\Lan_J \colon C^* \to D^*.
\end{equation}

By \cite[\S 4 eq. (4.19)]{kelly:enriched}, for any $G \in (D^\op)^*$ there is an isomorphism
\begin{equation}\label{coend-Lan}
G \otimes_D \Lan_JF = (J^*G) \otimes_C F.
\end{equation}
Given an object $c \in C$, we write $y_c$ for the functor $\map_C(c,-) \in C^*$.
It follows from \eqref{lkanadjoint} and Yoneda's Lemma \cite[2.4]{kelly:enriched} that
\begin{equation}\label{lancorep}
\Lan_J y_c = y_{Jc} = \map_D(Jc,-).
\end{equation}

Let $C_1,\ldots,C_n$ be $\Top$-categories. 
The \emph{external product} of $F_i \in C_i^*$ is
\begin{equation}\label{def-exttimes}
\bigexttimes_i F_i \in \Bigl(\prod_i C_i\Bigr)^*, \qquad \text{ defined by } \qquad (c_1,\ldots,c_n) \mapsto \prod_i F_i(c_i)
\end{equation}
where the product on the right is taken in $\Top$.
This gives a continuous functor $\bigexttimes \colon \prod_i C_i^* \to (\prod_i C_i)^*$. Given functors $J_i \colon C_i \to D_i$ between $\Top$-categories ($i=1,\dots,n$), we claim that
\begin{equation}\label{Lkan-product-functors}
\Lan_{\prod_i J_i}\Bigl(\bigexttimes_i F_i\Bigr) = \bigexttimes_i (\Lan_{J_i} F_i).
\end{equation}
Indeed, for any $(d_1,\ldots,d_n) \in \prod_i D_i$, it follows from Fubini's theorem and Yoneda's Lemma \cite[\S 3.3 and \S 3.10 eqs. (3.63), (3.67) , \S 2.4]{kelly:enriched} that
\begin{multline*}
\Lan_{\prod_i J_i}(\bigexttimes_i F_i)(d_1,\ldots,d_n) =
\int\limits^{(c_1,\ldots,c_n) \in \prod_i C_i} \big(\prod_i \map_{D_i}(J(c_i),d_i)\big) \times \prod_i F_i(c_i) 
\\
= \prod_i \Big( \int\limits^{c_i \in C_i} \map_{D_i}(J_i(c_i),d_i) \times F_i(c_i)\Big) = 
   \big(\bigexttimes_i \Lan_{J_i}F_i \big)(d_1,\ldots,d_n).
\end{multline*}
For any $c_i \in C_i$, the external product $\bigexttimes_i y_{c_i}$ is the functor 
$y_{(c_1,\ldots,c_n)} \in (\prod_i C_i)^*$.
Hence \eqref{lancorep} implies that for any functor $Q \colon \prod_i C_i \to D$,
\begin{equation}\label{lan-exttimes-corep}
\Lan_Q(\bigexttimes y_{c_i}) = y_{Q(c_1,\ldots,c_n)} = \map_D( Q(c_1,\ldots,c_n),-).
\end{equation}

Following Hollender and Vogt \cite{hollender-vogt}, we define the two-sided bar construction $\B(G,C,F)\in s\Top$ of 
functors $F\in C^*$ and $G \in (C^\op)^*$ as the simplicial space whose space of $n$-simplices is
\[
\B_n(G,C,F) = \coprod_{c_n,\dots,c_0 \in C} \biggl( G(c_0) \times \prod_{i=0}^{n-1} \map_C(c_{i+1},c_i)  \times F(c_n)\biggr).
\]
The \emph{nerve} of $C$ is the simplicial space $\B C=\B(*,C,*)$.
If $G\in (C^\op \times D)^*$ is a bifunctor, then $\B(G,C,F)$ inherits the structure of a functor $D \to s\Top$, 
and similarly for $F$. 

We would like to take the geometric realization of these bar constructions. In order for this to be homotopy meaningful, we need the following result:
\begin{lem}\label{reedycofibrantB}
If the values of $F$ and $G$ are CW-complexes then $\B(G,C,F)$ is cofibrant in the Reedy model structure on simplicial spaces (cf.~\cite[VII.\S 2]{goerss-jardine}). 
\end{lem}
\begin{proof}
We need to see that the inclusion of the latching objects
\[
\phi_n\colon L_n\B(G,C,F) \to \B_n(G,C,F)
\]
are cofibrations of topological spaces. We will in fact show that $\phi_n$ is the inclusion of a CW-subcomplex. Its image is given by
\[
\im(\phi_n) = 
\{ (g,f_n,\dots,f_1,f) \in \B_n(G,C,F)\mid \text{ at least one $f_i$ is an identity}\},
\] 
and the inclusion $\im (\phi_n) \subseteq \B_n(G,C,F)$ is an inclusion of a CW-subcomplex because we assume 
that all the spaces $\map_C(-,-)$ and the spaces $F(-)$ and $G(-)$ are CW-complexes and because $\{ \id_c \}$ is a $0$-cell in $\map_C(c,c)$. 
Thus, $\im(\phi_n) \subseteq \B_n(G,C,F)$ is a cofibration.

In order to prove that $\phi_n$ is a cofibration, it only remains to prove that it is a homeomorphism onto its image.
First, $\phi_n$ is injective because the simplicial set underlying $\B(G,C,F)$, as any simplicial set, is cofibrant. Furthermore, $\phi_n$ is a closed map: every degeneracy map $s_i$ is a closed inclusion since 
$\{ \id_c \} \subseteq \map_C(c,c)$ are closed inclusions, and $\phi_n$ is a quotient of the closed map 
\[
\coprod_{i=0}^{n-1} \B_{n-1}(G,C,F) \xto{\sum s_i} \B_n(G,C,F)
\]
by the definition of $L_n\B(G,C,F)$.
\end{proof}

Following \cite{hollender-vogt} we define the \emph{homotopy coend} 
to be the geometric realization $B(G,C,F) = |\B(G,C,F)|$. That is,
\[
B(G,C,F) = \Delta^\bullet \otimes_\Delta \B(G,C,F) = \int^{n \in \De} \Delta^n \times \B_n(G,C,F),
\]
where $\Delta^\bullet\colon \Delta \to \Top$ is the standard $n$-simplex functor.

It follows from \cite[Proposition VII.3.6]{goerss-jardine} that if
$F \to F'$ and $G \to G'$ are (objectwise) weak equivalences of functors 
then the resulting $B(G,C,F) \to B(G',C,F')$ is a weak equivalence of  cofibrant spaces.

Define
$$
\E C = \B(\map_C(-,-), C, *) \in s(C^*) \quad \text{and} \quad EC=|\E C| \in C^*.
$$
The Reedy cofibrant simplicial space $\E C(c)$ is augmented by $\E_{-1}C(c)=*$, and the inclusion of 
identity morphisms into $\map_C(c,c)$ gives rise to maps $s_{-1} \colon \E_n C(c) \to \E_{n+1} C(c)$ for all $n \geq -1$ by
\[
\map_C(c_0,c)\! \times \! \prod_{i=0}^{n-1} \! \map_C(c_{i+1},c_i) 
\to
   \map_C(c,c)\! \times \! \map_C(c_0,c) \! \times \! \prod_{i=0}^{n-1} \! \map_C(c_{i+1},c_i).
\]
The maps $s_{-1}$ form a left contraction for $\E C (c) \to *$, see \cite[\S 3.2 and Propositions 3.3 and 3.5]{cdi:02}, which implies that
$EC(c)$ is contractible because the geometric realization of the Reedy cofibrant simplicial space $\E C(c)$ is equivalent
to its homotopy colimit as a functor $\De^\op \to \Top$.

We remark that $\E C^{\op} \otimes_{C} F = \B(*,C,F)$ \cite[\S 3.10 eq. (3.72)]{kelly:enriched}.

\begin{defn}\label{holimdef}
The enriched homotopy (co)limit of $F \in C^*$ \cite{dror:diagrams} is defined by
\begin{eqnarray*}
\hholim{C} F &=& [EC,F]_C = \int_{c \in C} \hom(EC(c),F(c)),
\\
\hocolim_C F &=& (EC^{\op}) \otimes_C F = \int^{c \in C} EC^{\op}(c) \times F(c).
\end{eqnarray*}
\end{defn}
Here $\hom(-,-)$ is the mapping space of two objects in $\Top$.
Therefore
$$
\hocolim_{C} F = EC^{\op} \otimes_C F = |\B(*,C,F)|
$$
This construction is homotopy invariant as detailed in the next result.

\begin{prop}\label{weqhocolim}
Let $C,D$ be CW-categories.
\begin{enumerate}
\item
Let $X,Y \in C^*$ be functors whose values are CW-complexes.
A natural transformation $T \colon X \to Y$ induces a natural map $\hocolim_C X \to \hocolim_C Y$ which is
a weak equivalence if $T(c)$ is a weak equivalence $X(c) \to Y(c)$ for all $c \in C$.

\item
Fix a functor $X \in C^*$ whose values are CW-complexes and consider a continuous functor $F \colon D \to C$ between two categories $C$ and $D$ with the same objects. Assume that $F$ is the identity on objects and
induces weak equivalences on mapping spaces.
Then the natural map $\hocolim_D F^*X \to \hocolim_C X$ is a weak equivalence.
\end{enumerate}
\end{prop}

\begin{proof}
The natural transformation $T$ and the functor $F$ give rise to morphisms of Reedy cofibrant simplicial spaces
\[
\B(*,C,X) \xto{T} \B(*,D,Y) \qquad \text{ and } \qquad  
\B(*,D,F^*(X)) \xto{F} \B(*,C,X)
\]
which are homotopy equivalences in every simplicial degree by the hypotheses on $T$ and $F$.
We obtain weak equivalences by taking geometric realizations by \cite[Proposition VII.3.6]{goerss-jardine}.
\end{proof}

Similarly to Lemma~\ref{reedycofibrantB}, for any $F \in C^*$ the cosimplicial space $[\E C, F]_C$ is Reedy fibrant.
Here there is no restriction on the spaces $F(c)$ because every $X \in \Top$ is fibrant.
It follows that $\holim_C F= \Tot([\E C, F]_C)$ is homotopy invariant.

\begin{prop}\label{ECcofibrant}
The inclusion $\E_0 C \to |\E C| = EC$ induces a fibration in $\Top$
\[
\map_{C^*}(EC,X) \to \map_{C^*}(\E_0C,X).
\]
\end{prop}

For $c \in C$, let $y_c = \map_C(c,-) \in C^*$ be the continuous functor corepresented by $c$.
Given a space $A$ we write $y_c \otimes A$ for the functor $c' \mapsto y_c(c') \times A$.

\begin{lem}\label{freecofiblemma}
Let $A \to B$ be a cofibration and let $X \in C^*$ be a continuous functor. 
Then
\[
\map_{C^*}(B \otimes y_c,X) \to \map_{C^*}(A \otimes y_c,X)
\]
is a fibration for every $c \in C$.
\end{lem}
\begin{proof}
By adjunction and Yoneda's Lemma \cite[\S 3.10 eq. (3.72)]{kelly:enriched}, the map under consideration is isomorphic to 
$\map(B,X(c)) \to \map(A,X(c))$, which is a fibration by hypothesis on $A \subseteq B$. One can check this
by induction on the skeleta of the relative CW-complexes $(B,A)$ and applying \cite[\S 2.8, Theorem 2]{spanier:at} 
to the inclusions $\partial \De^n \subseteq \De^n$.
\end{proof}

\begin{proof}[Proof of Proposition \ref{ECcofibrant}]
Let $\sk_n\E C$ denote the (objectwise) $n$-skeleton of the simplicial $C$-space $\E C$. 
By \cite[VII.3.8]{goerss-jardine} and Lemma~\ref{reedycofibrantB} there are pushout squares 
\begin{equation}\label{freecofib:1}
\begin{matrix}
\xymatrix@1{
\coprod\limits_{\mathbf{c}} \left(A_{\mathbf{c}} \times \Delta^n \smash{\underset{A_{\mathbf{c}}\times \partial\Delta^n}\cup}  B_{\mathbf{c}} \times \partial\Delta^n \right) \otimes y_{c_0} \ar[r] \ar[d] & |\sk_{n-1}\E C| \ar[d]\\
\coprod\limits_{\mathbf{c}} B_{\mathbf{c}} \times \Delta^n \otimes y_{c_0} \ar[r] & |\sk_n\E C|,
}
\end{matrix}
\end{equation}
where $\mathbf{c} = (c_n,\dots,c_0)$, $B_{\mathbf{c}} = \prod_{i=0}^{n-1} \map_C(c_{i+1},c_i)$, and 
$A_{\mathbf{c}} \subseteq B_{\mathbf{c}}$ is the subset of degenerate $n$-simplices, 
i.e. those where at least one map is the identity. 
We have seen in Lemma~\ref{reedycofibrantB} that 
$\B C$ is Reedy cofibrant and, in fact, that the map $L_n\B C \to \B_n C$ is an inclusion of CW-complexes.
Moreover, $L_n\B C = \coprod A_{\mathbf{c}}$ and $\B_n C = \coprod_{\mathbf{c}} B_{\mathbf{c}}$.
This implies that  the left hand side of \eqref{freecofib:1} satisfies the conditions of Lemma~\ref{freecofiblemma}.
Thus, by applying $\map_C(-,X)$ to \eqref{freecofib:1}, it becomes a pullback square 
in which one of the sides is a fibration and therefore the side opposite to it, namely 
\[
\map_{C^*}(|\sk_n \E C|,X) \to \map_{C^*}(|\sk_{n-1} \E C|,X)
\]
is also a fibration.
Since the inverse limit of a tower of fibrations is a fibration and since 
$EC = \colim_n |\sk_n \E C|$, we obtain the required  
fibration $\map_{C^*}(EC,X) \to \map_{C^*}(\sk_0 \E C,X) = \map_{C^*}(\E_0C,X)$.
\end{proof}

\section{\texorpdfstring{$A_\infty$}{A-infinity}-algebras, \texorpdfstring{$A_\infty$}{A-infinity}-monads, and their completions}
\label{sec-a-monads-completion}

\subsection{Operads and associated categories}

\begin{defn}
By an \textbf{operad} $A$ of topological spaces we always mean a non-symmetric operad as in \cite[Definition 3.12]{may:72}.
That is, $A$ consists of a sequence of spaces $A(n)$ for $n \geq 0$ with associative composition operations
\[
A(n) \times A(k_1) \times \dots \times A(k_n) \to A\Bigl(\sum k_i\Bigr)
\]
and a base point $\iota \in A(1)$ which serves as a unit.
We call $A$ an \emph{$A_\infty$-operad} if all the spaces $A(n)$ are contractible CW-complexes. 
\end{defn}

\begin{nota}\label{operadnotation}
The following coordinate-free description of an operad $A$ will be useful: 
If $S$ is a finite ordered set, we write $A(S)$ for $A(\#S)$. 
If $\vp\colon S \to T$ is a monotonic map of finite ordered sets, we write
\[
A(\vp) = \prod_{t \in T} A(\vp^{-1}t), \quad \text{where $\vp^{-1}t = \vp^{-1}(\{t\}) \subseteq S$}.
\]
Thus, the operad structure is given by maps $\mu(T,\vp)\colon A(T) \times A(\vp) \to A(S)$, or equivalently, 
by a collection of maps
\[
\mu(\psi,\vp)\colon A(\psi) \times A(\vp) \to A(\psi \circ \vp) \quad \text{for all } S \xto{\vp} T \xto{\psi} W
\]
which are given by 
\begin{align*}
\prod_{w \in W} A(\psi^{-1}w) \times \prod_{t \in T} A(\vp^{-1}t) =& \prod_{w \in W} \biggl(A(\psi^{-1}w) \times A(\vp|_{\vp^{-1}\psi^{-1}w})\biggr)\\
\xto{\prod \mu(\psi^{-1}w,\vp|_{\vp^{-1}\psi^{-1}w})} & \prod_{w \in W} A(\vp^{-1}\psi^{-1}w). 
\end{align*}
and which are associative and unital.
The base points in $A(\{s\})$ assemble to a base point $\iota_S \in A(\id_S)$ which acts as a unit for $\mu$.
\end{nota}

Let $\Deplus$ be the essentially small category of finite ordered sets and monotonic maps, 
where we allow the empty set as an object, and $\De$ the full subcategory of nonempty finite ordered sets. 
We might as well take $\De$ and $\Deplus$ to be the small skeleton of ordered sets $[n] = \{0,\dots,n\}$, 
but the coordinate-free setting is more natural.
Concatenation of ordered sets makes $\Deplus$ into a (non-symmetric) monoidal category; we will denote this monoidal structure by $\sqcup$, keeping in mind that it is not a categorical coproduct. 

We will now ``thicken up'' the categories $\De$ and $\Deplus$ 
by allowing morphisms to be parameterized by the spaces of an operad. Such categories were called ``categories of operators in standard form'' in \cite{boardman-vogt:homotopy-everything}.

\begin{defn}\label{defdeltaA}
For an operad $A$, let $\De(A)$ be its associated category of operators \cite{boardman-vogt:homotopy-everything,may-thomason} and let $\Deplus(A)$ be the obvious enlarged category containing the empty set. Explicitly,
\[
\map_{\Deplus(A)}(S,T) = \coprod_{\vp \in \Hom_{\Deplus}(S,T)} A(\vp) \qquad \text{(cf. Notation~\ref{operadnotation})}
\]
\end{defn}
In \cite[\S4]{may-thomason} it was shown that this is indeed a monoidal topological category over the category $\Deplus$ of finite ordered set with concatenation $\sqcup$ of ordered sets. Namely, the maps $A(S) \to *$ give rise to a functor
\[
\pi \colon \Deplus(A) \to \Deplus.
\]

If $A$ is $A_\infty$ then $\Deplus(A)$ is a CW-category (\ref{toplogicalcategories}). The following proposition is immediate from the definitions:

\begin{prop} \label{operadequivalence}
\begin{enumerate}
\item 
$\map_{\Deplus(A)}(\emptyset,S)=\prod_S A(\emptyset)$ for any $S \in \Deplus(A)$.
If $S \neq \emptyset$ then $\map_{\Deplus(A)}(S,\emptyset)=\emptyset$.

\item $\Deplus(\assoc) = (\Deplus,\sqcup)$ where $\assoc$ denotes the ``associative operad'', 
namely the operad all of whose spaces are singletons. \label{assocdeplus}

\item If $A$ is an $A_\infty$-operad then $\pi \colon \Deplus(A) \to \Deplus$ induces homotopy equivalences on the mapping
spaces.\qed \label{ainftyequiv}
\end{enumerate}
\end{prop}

\subsection{Algebras over operads in monoidal categories} \label{endcblog}

Let $(\E,\comp,I)$ be a monoidal category enriched over $\Top$ such that $\comp$ is an enriched functor.
For an object $K$ we write $K^{\comp n}$ for $K \comp \cdots \comp K$ ($n$ times) where $K^{\comp 0}$ is the monoidal identity object $I$.
\begin{defn}\label{oper-of-obj}
The object $K$ gives rise to an operad $\opEND(K)$ defined as follows. For a totally ordered set $S$, the $S$th space $\opEND(K)(S)$ is $\map_\E(K^{\comp S},K)$. For a monotonic map $\vp\colon S \to T$, the composition law is given by (cf. \cite[Definition~1.2]{may:72})
\begin{multline*}
\map_\E(K^{\comp T},K) \times \prod_{t \in T} \map_{\E}(K^{\comp \vp^{-1}t},K) \xto{\text{monoidal}} 
\\
\map_\E(K^{\comp T},K) \times \map_\E(K^{\comp S},K^{\comp T}) 
\xto{\text{compose}} \map_\E(K^{\comp S},K).
\end{multline*}
If $A$ is an operad, an \emph{$A$-algebra in $(\E,\comp,I)$} is an object $K \in \E$ with a morphism of operads $A \to \opEND(K)$.
\end{defn}

\begin{remark}
In the classical case of $A$-algebras in $(\Top,\times,*)$, there is an associated monad $\mathbb A$ such that $\mathbb A$-algebras are the same as $A$-algebras in the above sense. \emph{This does not hold in general monoidal categories.} In fact, there is no such thing as a free $A$-algebra on an object $K$. The reason for this is the failure of the monoidal structure to be ``linear'' over spaces: Both of the following conditions are satisfied in $\Top$ but not in a general monoidal category (the second property assumes that the category is tensored over $\Top$):
\begin{itemize}
\item $K \comp (L_1 \sqcup L_2) \cong (K \comp L_1) \sqcup (K \comp L_2)$
\item $K \comp (X \otimes L) \cong X \otimes (K \comp L)$ for all spaces $X$.
\end{itemize}
\end{remark}

\begin{defn}\label{def-k-completion}
Let $K$ be an $A$-algebra in $(\E,\comp,I)$ with structure map $\si \colon A \to \opEND(K)$.
The $A$-\emph{cobar construction} of $K$ is the monoidal functor
\[
R_K \colon (\Deplus(A),\sqcup,\emptyset) \to (\E,\comp,I)
\]
defined by setting $R_K(S)= K^{\comp S}$, and on morphism spaces (cf. Notation \ref{operadnotation}) by
\begin{align*}
\map_{\Deplus(A)}(S,T) = &\coprod_{\vp \in \Deplus(S,T)} \prod_{t \in T} A(\vp^{-1}t)\\
\xto{\coprod_\vp\prod_t \si(\vp^{-1}t)} &
\coprod_{\vp \in \Deplus(S,T)} \prod_{t \in T} \map(K^{\comp \vp^{-1}t},K) 
\xto{\ \comp \ }
\map(K^{\comp S},K^{\comp T}).
\end{align*}
\end{defn}

The continuity of the maps between the mapping spaces follows from the continuity of the operad map $\si$ 
and the continuity of the monoidal structure maps in $\E$.
It clearly carries identities in $\Deplus(A)$ to identities in $\E$ because $\si$ carries the
base point of $A(\{s\})$ to the base point of $\map(K,K)$.
To see that $R_K$ respects composition, let us consider morphisms $S \xto{a} T$ and $T \xto{b} W$ in $\Deplus(A)$.
Thus, $a=(a_t)_{t \in T}$ and $b=(b_w)_{w \in W}$ where $a_t \in A(\vp^{-1}t)$ and $b_w \in A(\psi^{-1}w)$ for some monotonic functions $\vp \colon S \to T$ and $\psi \colon T \to W$.
Denote $f_t = \si(a_t) \in \map(K^{\comp \vp^{-1}t},K)$ and $g_w = \si(b_w) \in \map(K^{\comp \psi^{-1}w},K)$.
Let $\bullet$ denote the composition operation in $A$.
Since $\si$ is a morphism of operads it follows from Definition~\ref{oper-of-obj} of $\opEND$, Definition~\ref{def-k-completion}, and by writing $T=\bigsqcup_w \psi^{-1}w$ that
\begin{multline*}
R_K(b \circ a) = R_K\bigl((b_w \bullet (a_t)_{t \in \psi^{-1}w})_{w \in W}\bigr) = 
\bigcomp_{w \in W} \bigl( \si(b_w \bullet (a_t)_{t \in \psi^{-1}w}) \bigr) 
\\
= \bigcomp_{w \in W} \biggl(g_w \circ \Bigl(\bigcomp_{t \in \psi^{-1}w} f_t\Bigr)\biggr) =
\Bigl(\bigcomp_{w \in W} g_w\Bigr) \circ \Bigl(\bigcomp_{t \in T} f_t\Bigr) = R_K(b) \circ R_K(a).
\end{multline*}
The fact that the monoidal operation $\comp$ in $\E$ is associative is the reason that $R_K$ 
is a monoidal functor.
The details are straightforward and are left to the reader.

\begin{example}\label{examplemonads}
An ordinary monoid $K$ in $\E$ is an $\assoc$-algebra. Indeed the images of $\assoc_0$ and $\assoc_2$ in $\opEND(K)$ give the unit $\eta\colon I \to K$ and the monoid operation
$\mu\colon K \circ K \to K$.
Inspection reveals that $R_K \colon \Deplus \to \E$ is the standard cobar construction with $(R_K)^n = K^{\comp(n+1)}$ and with coface maps $d^i=K^i \eta K^{n-i}$ and $s^i = K^i \mu K^{n-i}$. 
\end{example}

\begin{defn}
\label{def-completion}
Given an operad $A$, let $J \colon \De(A) \to \Deplus(A)$  denote the inclusion of the full subcategory $\De(A)$
spanned by the non-empty sets.
The \emph{completion} $\hat K\in \E$ of an $A$-algebra $K$ is 
\[
\hat K = \hholim{\De(A)} J^*R_K = [E\De(A),J^*R_K]_{\De(A)} \quad (\text{cf. Definition~\ref{holimdef}.})
\]
\end{defn}

\subsection{Endofunctor categories and \texorpdfstring{$A$-monads}{A-monad}}\label{ss:endofunctorCat}

Our main examples for $\E$ are categories of endofunctors. Let $\C$ be a category which is complete and cocomplete, and which is tensored by
$\otimes^\C \colon \C \times \Top \to \C$ and cotensored by $[-,-]^\C \colon \Top^\op \times \C \to \C$.
The adjoint $\map_\C \colon \C^\op \times \C \to \Top$ of $\otimes^\C$ endows $\C$ with
an enrichment over $\Top$.

Let $\E=\END(\C)$ denote the Category of all the continuous endofunctors of $\C$ together with natural transformations between them.
The word Category is capitalized because the natural transformations between two endofunctors need not form a set. To overcome this difficulty we pass to a larger universe \cite[I.6]{maclane:cwm}.
The morphism Spaces in $\END(\C)$ are denoted $\Map(K,L)$.

\begin{remark}
Throughout this paper we will not make an essential use of the larger universe.
The introduction of the Category $\END(\C)$ is only done for the sake of book-keeping and we will never refer to morphism Sets in $\END(\C)$.
One may prefer to consider $\END(\C)$ as a convenient notation for the classes of continuous endofunctors $K$ on $\C$ and the classes of natural transformations between them.
In fact, we are only going to construct explicit functors $D \to \END(\C)$ from small topological categories $D$. 
These are the same as continuous (honest) functors $F \colon D \times \C \to \C$.
Alternatively, these are assignments of endofunctors $F(d)$ for every object of $D$ and of natural transformations $F(d) \to F(d')$ for every morphism of $D$ such that the compositions
$$
\map_D(d,d') \to \Nat(F(d),F(d')) \xto{\ev_X} \map_\C(F(d)(X),F(d')(X))
$$
are continuous functions for all objects $X \in \C$.
It now becomes obvious that it is just a matter of convenience to assume that $\END(\C)$ is a Category (capitalized or not).
\end{remark}

The Category $\END(\C)$ is tensored and cotensored over $\Top$, where for every space 
$A \in \Top$ and every $F \in \END(\C)$ we define $F \otimes A$ and $[A,F]$ in $\END(\C)$ as the endofunctors defined by $x \mapsto F(x) \otimes^\C A$ and $x \mapsto [A,F(x)]^\C$.
The adjoint of the tensor is $\Map(F,G) = \int_{X \in C} \map_\C(F(X),G(X))$.

Composition of functors equips $\END(\C)$ with an enriched monoidal structure $(\comp,I)$:
\[ 
K \comp L = K \circ L \quad \text{and} \quad I = \id_\C.
\]
We will reserve the symbols $\circ$ and $\id$ for composition of natural transformations (i.e. Morphisms) and the identity natural transformation in this Category.
Explicitly, 
$$
\Map(K,L) \times \Map(K',L') \xto{\comp} \Map(K \comp K', L \comp L') 
$$ 
has the effect 
$(f,f') \mapsto f \comp f' = (L \comp f') \circ (f \comp K') = (f \comp L') \circ (K \comp f')$. The continuity of $K$ and $L$ guarantees the continuity of the map above between the mapping Spaces.
To see that this gives rise to a monoidal structure one calculates that 
$\id_K \comp \id_L = (L \comp \id_K) \circ (\id_L \comp K) = \id_{K \comp L} \circ \id_{K \comp L}=\id_{K \comp L}$.
Also, given $f_i \colon K_i \to L_i$ and $g_i \colon J_i \to K_i$ where $i=1,2$ we have
\begin{multline*}
(f_1 \comp f_2) \circ (g_1 \comp g_2) = (f_1 \comp L_2) \circ (K_1 \comp f_2) \circ (K_1 \comp g_2) \circ (g_1 \comp J_2) \\ 
=
(f_1 \comp L_2) \circ (K_1 \comp (f_2 \circ g_2)) \circ (g_1 \comp J_2) 
\\
=(f_1 \comp L_2) \circ (g_1 \comp L_2) \circ (J_1 \comp (f_2 \circ g_2)) = (f_1 \circ g_1) \comp (f_2 \circ g_2).
\end{multline*}

Any $K \in \END(\C)$ gives rise to an Operad $\opEND(K)$ as in Definition~\ref{oper-of-obj},
whose $n$th Space is $\Map(K^{\comp n},K)$.

\begin{defn}\label{def-a-monad}
An \emph{$A$-monad} in $\C$ is an $A$-algebra in $(\END(\C),\comp,I)$.
\end{defn}

Thus, an $A$-monad in $\C$ is an endofunctor of $\C$, namely an object $K \in \END(\C)$ together with a Morphism of Operads $\si \colon A \to \opEND(K)$, see Definition~\ref{def-k-completion}, where $\opEND(K)$ is an Operad of Spaces and $A$ is an operad of (small) spaces.
If one wishes to avoid the set-theoretic difficulties around $\END(\C)$ mentioned before, one can replace the map of Operads $\si\colon A \to \opEND(K)$ with a collection of natural transformations of endofunctors of $\C$
\[
A(n) \otimes^\C K^{\comp n} (-) \xto{ \al(n) } K(-)
\]
which make $K$ an algebra over the operad $A$ in a sense similar to \cite{may:72}.
To do this, one needs to use the continuity of $K$, which gives rise to natural maps $U \otimes^\C K(-) \to K(U \otimes^\C -)$ for any topological space $U$.

Note that ordinary monads \cite{eilenberg-moore:65} (called triples there) are $\assoc$-monads in our terminology, and the cobar construction $R_K(X)$ of Definition~\ref{def-k-completion} and completion $\hat K(X)$ of Definition~\ref{def-completion} coincide with the usual definitions \cite{yellowmonster}.

\section{The extended \texorpdfstring{$A_\infty$}{A-infinity}-operad}
\label{ahp-section}

The goal of this section is to construct an ``extended'' operad $\AHP$ associated with any operad $A$. In Section~\ref{action-on-completion}, this operad will be seen to act naturally on the completion of any $A$-algebra (Definition~\ref{def-completion}.)

Throughout this section we fix an operad $A$ and let $\Dplus$ denote its category of operators $\Deplus(A)$ (Definition~\ref{defdeltaA}).
We denote by $J \colon D \to \Dplus$ the inclusion of the full subcategory spanned by the objects 
$\emptyset \neq S \in \Dplus$.
The $S$-fold product of a category $C$ with itself, where $S$ is a finite ordered set, is denoted $C^S$.

\begin{defn}\label{ord-sum-functors}
For a finite ordered set $S$, the ordinal sum functor is defined by
$$
u(S) \colon {\Dplus}^S \to \Dplus, 
\qquad \qquad
(V_s)_{s \in S} \mapsto \bigsqcup_{s \in S} V_s.
$$
\end{defn}

The monoidal axioms for $\sqcup$ readily imply
\begin{lem}\label{ordsumeasy}
\begin{enumerate}
\item $u(\emptyset)\colon  * \to \Dplus$ is the inclusion of the $\sqcup$-unit $\emptyset$ into $\Dplus$. \label{uemptyset}

\item $u(*) = \id_{\Dplus}$.  \label{uonepoint}

\item $u(S) = u(T) \circ \Bigl(\underset{t \in T}{\prod} u(\vp^{-1}t)\Bigr)$ for any monotonic $S \xto{\vp} T$ (see Notation \ref{operadnotation}).  \label{uS-in-steps} \qed
\end{enumerate}
\end{lem}

The functor $E\D \in \D^*$ (see Section~\ref{s:enriched-categories}) can be extended to a continuous functor $\TED\in \Dplus^*$ by setting $\TED(\emptyset) = \emptyset$.
In fact, by Proposition \ref{operadequivalence} and  \cite[Proposition 4.23]{kelly:enriched}
\begin{equation}
\label{TEDasLan}
\TED=\Lan_J E\D
\end{equation}
because $J$ is fully faithful.

\begin{defn}
\label{def-epsilonS}
For a finite ordered set $S$ define $\epsilon(S)\in \Dplus^*$ by
\[
\epsilon(S) = \Lan_{u(S)} \Bigl( \bigexttimes_{s \in S}\TED\Bigr) \qquad \text{(see Definition~\ref{def-exttimes}).}
\]
\end{defn}

Let $\vp \colon S \to T$ be a monotonic map and write $S= \bigsqcup_{t \in T} \vp^{-1}t$.
It follows from \eqref{Lkan-product-functors} and from \cite[Theorem~4.47]{kelly:enriched} that
\begin{equation}\label{epsilon-of-sqcup}
\begin{aligned}\epsilon(S) = &\Lan_{u(S)}(\bigexttimes_S \TED) \underset{\text{Lemma \ref{ordsumeasy}}}= 
\Lan_{u(T)} \circ \Lan_{\prod_t u(\vp^{-1}t)} \Bigl(\bigexttimes_{t \in T} \bigexttimes_{\vp^{-1}t} \TED\Bigr) \\
\underset{\text{\eqref{Lkan-product-functors}}}= &\Lan_{u(T)} \Bigl( \bigexttimes_t \Lan_{u(\vp^{-1}t)}(\bigexttimes_{\vp^{-1}t}\TED) \Bigr) = 
\Lan_{u(T)} \Bigl( \bigexttimes_{t \in T} \epsilon(\vp^{-1}t) \Bigr)
\end{aligned}
\end{equation}

\begin{defn}\label{construction-AHP}
Given a finite ordered set $S$, define a space $\AHP(S)$ by
\[
\AHP(S) = \map_{\Dplus^*}(\TED,\epsilon(S)).
\]
Analogously to Notation~\ref{operadnotation}, for a monotonic map $\vp \colon S \to T$ we define $\AHP(\vp)$ by
\[
\AHP(\vp) = \prod_{t \in T} \AHP(\vp^{-1}t) = \prod_{t \in T} \map_{\Dplus^*}(\TED,\epsilon(\vp^{-1}t))
\]
If $|S|=1$ then Lemma~\ref{ordsumeasy}(\ref{uonepoint}) implies that $\epsilon(S)=\TED$. 
The identity on $\TED$ equips $\AHP(S)$ with a natural base point.
For any $\vp \colon S \to T$, define a continuous function 
\begin{equation}
\label{Gamma-T-phi}
\tilde\mu(T,\vp) \colon \AHP(T) \times \AHP(\vp) \to \AHP(S)
\end{equation}
as follows.
For any element $f \in \AHP(T)$ and any element $(g_t)_{t \in T} \in \AHP(\vp)$ 
recall that $f$ and $g_t$ are
morphisms $f \colon \TED \to \epsilon(T)$ and $g_t \colon \TED \to \epsilon(\vp^{-1}t)$.
Apply \eqref{epsilon-of-sqcup} to define
\[
\tilde{\mu}(T,\vp) \colon (f,(g_t)) \mapsto \Lan_{u(T)}(\bigexttimes_{t \in T} g_t) \circ f.
\]
The continuity of $\tilde{\mu}(T,\vp)$ follows from the fact that $\Lan_{u(T)}$ and $\bigexttimes_T$ are 
continuous functors (\ref{lancontinuous}, \ref{def-exttimes}).

Following Notation~\ref{operadnotation}, this gives rise to maps 
\begin{equation}
\label{def-Gamma-psi-phi}
\AHP(\psi) \times \AHP(\vp) \xto{\tilde\mu(\psi,\vp)} \AHP(\psi \circ \vp) \qquad 
\text{for all $S \xto{\vp} T \xto{\psi} W$.}
\end{equation}
That is, $\tilde{\mu}(\psi,\vp)=(\tilde{\mu}(\psi^{-1}w,\vp|_{\vp^{-1}\psi^{-1}w}))_{w \in W}$
\end{defn}

The goal of this section is to prove the following two results.

\begin{thm}\label{ahpoperad}
The spaces $\AHP(S)$ together with the maps $\tilde{\mu}(\vp,\psi)$ give rise to an operad $\AHP$ which we call
the ``extended operad'' associated to $A$.
\end{thm}

\begin{thm}\label{ahpcontractible}
If $A$ is an $A_\infty$-operad then the spaces $\AHP(S)$ are weakly contractible.
\end{thm}

\begin{proof}[Proof of Theorem~\ref{ahpoperad}]
Given $\vp \colon S \to T$ we denote an element of $\AHP(\vp)$ by $\fff = (f_t)_{t \in T}$ for 
$f_t \in \map_{\Dplus^*}(\TED,\epsilon(\vp^{-1}t))$.

For monotonic maps $S \xto{\vp} T \xto{\psi} W$ and elements $\fff=(f_t) \in \AHP(\vp)$ and $\ggg=(g_w) \in \AHP(\psi)$, \eqref{Gamma-T-phi} and \eqref{def-Gamma-psi-phi} imply that
$\ggg \circ \fff := \tilde{\mu}(\ggg,\fff) = (h_w)$ where
\begin{equation}
\label{elementwise-composition}
h_w = \Lan_{u(\psi^{-1}w)} \Bigl(\bigexttimes_{t \in \psi^{-1}w} f_t\Bigr) \circ g_w\colon \TED \to \epsilon(\vp^{-1}\psi^{-1}w)
\end{equation}

\noindent \emph{Associativity of composition.} 
Let $S \xto{\vp} T \xto{\psi} W \xto{\theta} R$ be maps in $\Deplus$ and consider $\fff=(f_t) \in \AHP(\vp)$,
$\ggg=(g_w) \in \AHP(\psi)$, and $\hhh=(h_r) \in \AHP(\theta)$. 
Then
\begin{align*}
((\hhh\circ \ggg) \circ \fff)_r \underset{\eqref{elementwise-composition}}{ \; = \; } & \
   \Lan_{u(\psi^{-1}\theta^{-1}r)}\Bigl(\bigexttimes_{t\in\psi^{-1}\theta^{-1}r} f_t\Bigr) \circ 
   \biggl(\Lan_{u(\theta^{-1}r)}\Bigl(\bigexttimes_{w \in \theta^{-1}r} g_w\Bigr)\biggr) \circ h_r \\
\underset{\text{Lemma \ref{ordsumeasy}}}{\qquad = \qquad} & 
   \Lan_{u(\theta^{-1}r)} \circ \!\Lan_{\underset{{w\in\theta^{-1}r}}{\prod} \!\!\! u(\psi^{-1}w)} \Bigl(\bigexttimes_{w \in \theta^{-1}r} (\bigexttimes_{t \in \psi^{-1}w} f_t)\Bigr) \!\circ\! 
   \biggl(\Lan_{u(\theta^{-1}r)}\Bigl(\bigexttimes_{w \in \theta^{-1}r} g_w\Bigr)\biggr) \!\circ\! h_r \\
\underset{\text{\eqref{Lkan-product-functors}}}{\qquad=\qquad} & \Lan_{u(\theta^{-1}r)} \biggl( \bigexttimes_{w\in\theta^{-1}r}\Bigl(\Lan_{u(\psi^{-1}w)} (\bigexttimes_{t\in\psi^{-1}w} f_t)\Bigr) \circ g_w\biggr) \circ h_r\\
\underset{\eqref{elementwise-composition}}{ \qquad=\qquad} & 
\Lan_{u(\theta^{-1}r)}\Bigl( \bigexttimes_{w\in\theta^{-1}r}(\ggg \circ \fff)_w\Bigr) \circ h_r \quad =\quad (\hhh \circ(\ggg\circ \fff))_r.
\end{align*}

\noindent \emph{Unitality of composition.} There are two cases to check.
First consider $S \xto{\id_S} S \xto{\vp} T$ and fix $\fff=(f_t) \in \AHP(\vp)$.
Since $\Lan$ and $\bigexttimes$ are functors,
it follows from Definition~\ref{def-epsilonS} that
$$
(\fff \circ \iota_S)_t = \Lan_{u(\vp^{-1}t)} \Bigl(\bigexttimes_{s\in \vp^{-1}t} \id_{\TED}\Bigr) \circ f_t
= \id_{\epsilon(\vp^{-1}t)} \circ f_t = f_t. 
$$
Thus, $\fff \circ \iota_S= \fff$.
Next, consider $S \xto{\vp} T \xto{\id_T} T$ and and element $\fff =(f_t) \in \AHP(\vp)$.
By Lemma \ref{ordsumeasy},  $u(*):\Dplus\to\Dplus$ is the identity, hence
$$
(\iota_T \circ \fff)_t =  (\Lan_{u(*)} f_t) \circ \id_{\TED} = f_t.
$$
Therefore $\iota_T \circ \fff = \fff$.
\end{proof}

\begin{lem}\label{contractible-ordinalsum}
Fix $n \geq 1$ and let $\De^{\times n}$ be the $n$-fold product.
Then for any $T \in \De$ we have 
$$
\hhocolim{(P_1,\ldots,P_n) \in (\De^{\times n})^{\op} } \, \De(P_1 \sqcup \cdots \sqcup P_n,T) \ \simeq \ *.
$$
\end{lem}

\begin{proof}
When $n=1$ this is a triviality because 
$$
\hhocolim{\De^\op} \, \De(-,T) = \De(-,T) \otimes_{\De^\op} E\De = E\De(T) \simeq *. 
$$
Given $\th \in \De(P,T)$ we denote by $\max (\th) \in T$ the maximal element in the image of $\th$.
For any $t \in T$ we denote by $T-t$ the subset of $T$ consisting of all the elements $t' \in T$ such that $t'\geq t$.
We observe that
$$
\De(P \sqcup Q,T) =  \coprod_{\th \in \De(P,T)} \De(Q,T-\max(\th)).
$$
The lemma now follows by induction on $n$ because
\begin{align*}
\hhocolim{(P_1,\ldots,P_n) \in \De^{\times n}} \De\Big(\bigsqcup_{i=1}^n P_i,T\Bigr)
= & \hhocolim{P_1 \in \De}  \hhocolim{(P_2,\ldots,P_n)\in \De^{\times(n-1)}} \, \De\Bigl(P_1 \sqcup \bigsqcup_{i=2}^n P_i,T\Bigr) \\
=& \hhocolim{P_1 \in \De} \negthickspace \coprod_{\th \in \De(P_1,T)} 
\underbrace{\hhocolim{(P_2,\ldots,P_n)\in \De^{\times(n-1)}}  \De\Bigl(\bigsqcup_{i=2}^n P_i,T-\max(\th)\Bigr)}_{\simeq * \ \text{by induction}} \\
\simeq &
\hhocolim{P_1 \in \De} \De(P_1,T) \simeq *.
\end{align*}\qedhere
\end{proof}

\begin{proof}[Proof of Theorem~\ref{ahpcontractible}]
Since $\TED=\Lan_{J}ED$, equation \eqref{lkanadjoint} implies
$$
\AHP(S) = \map_{\Dplus^*}(\TED,\epsilon(S)) = \map_{D^*}(ED,J^*(\epsilon(S))) = \holim_{D} J^*(\epsilon(S)).
$$
It is therefore sufficient to prove that the values of $J^*(\epsilon(S))$ are contractible spaces, that is, 
we have to show that $\epsilon(S)(T) \simeq *$ for every $T \neq \emptyset$.

First, we note that $\bigexttimes_S ED = ED^S$ because $\E(D^S)=(\E D)^S$ and because geometric realization commutes with finite products by \cite[Theorem 11.5]{may:72}. From Definition \ref{def-epsilonS}, \eqref{Lkan-product-functors} and \cite[Theorem~4.47]{kelly:enriched} it follows that
$$
\epsilon(S) = \Lan_{u(S)}\Bigl(\bigexttimes_S \Lan_J ED\Bigr) = \Lan_{u(S)}\circ \Lan_{\prod_S J}\Bigl(\bigexttimes_S ED\Bigr) =
\Lan_{u(S) \circ (\prod_S J)} (ED^S).
$$
Set $F= u(S) \circ (\prod_S J)$, thus $\epsilon(S)=(\Lan_F ED^S)$.
By the definition of $\Lan$ and by Fubini's theorem in enriched categories \cite[\S 2.1 and 3.3]{kelly:enriched},
\[
\epsilon(S)(T) = \map_{\Dplus}(F(-),T) \otimes_{D^S} ED^S 
=  \hocolim_{(D^S)^\op} \,\map_{\Dplus}(F(-),T).
\]
There is an obvious ordinal sum functor $u'(S) \colon \Deplus^S \to \Deplus$ and inclusion $J' \colon \De \to \Deplus$ 
which render the following diagrams commutative:
$$
\xymatrix@1{
{\Dplus}^S \ar[d]_{\prod_S \pi} \ar[r]^{u(S)} &
\Dplus \ar[d]^{\pi} 
\\
{\Deplus}^S \ar[r]_{u'(S)} & \Deplus
}
\qquad \qquad
\xymatrix@1{
D \ar[r]^{J} \ar[d]_{\pi} &
\Dplus \ar[d]^{\pi}
\\
{\De} \ar[r]_{J'} & \Deplus.
}
$$
The vertical arrows induce the identity on object sets and are weak equivalences of categories by Proposition~\ref{operadequivalence}\eqref{ainftyequiv}.
Set $F':=u'(S) \circ (\prod_S J')$.
This is the functor $\De^{\times n} \to \Deplus$ where $(P_1,\ldots,P_n) \mapsto \bigsqcup_i P_i$. 
There results a weak equivalence of functors $\Dplus \to \Top$
\[
\map_{\Dplus}(F(-),T) \xto{\ \simeq \ }  (\prod_S \pi)^*\big(\Deplus(F'(-),T) \big).
\]
It follows from Proposition~\ref{weqhocolim} that there is a weak equivalence
\[
\epsilon(S)(T) \simeq 
\hhocolim{(\De^S)^\op} \Deplus(F'(-),T) 
\]
which is weakly equivalent to a point by Lemma \ref{contractible-ordinalsum} and the fact that $T \neq \emptyset$ and that
$\De$ is a full
subcategory of $\Deplus$.
\end{proof}

\section{The action of \texorpdfstring{$\AHP$}{A-tilde} on the completion} \label{action-on-completion}

As in Subsection~\ref{endcblog}, let $(\E,\comp,I)$ be a complete and cocomplete monoidal category compatibly tensored and cotensored over $\Top$.

Throughout this section we will write $f \colon X \to Y$ for the elements of 
$\map_\E(X,Y)$, cf. \cite[1.3]{kelly:enriched}. Furthermore, for $F,\; G \in \E$ and $A \in \Top$ we will write $F^A$ for $[A,F]$ and $\{ F,G \}$ for $\map_\E(F,G)$. 

In this section we will prove:

\begin{thm} \label{extendedoperadaction}
Let $A$ be an $A_\infty$-operad and $K$ be an $A$-algebra in $(\E,\comp,I)$ (cf. Definition~\ref{def-a-monad}). 
Then $\hat K$ is an $\AHP$-algebra where $\AHP$ is defined in Definition~\ref{construction-AHP}.
\end{thm}

Given $A,B \in \Top$ and  $F,G \in \E$ we consider the map $\la_{A,B}^{F,G}$ defined by 
\begin{equation}\label{defrhoadjoint}
A \times B \xto{\ev_A^{\#} \times \ev_B^{\#}}
\{F^A,F\} \times \{G^B,G\} \xto{\, \Diamond \, } \{F^A \comp G^B,F \comp G\}.
\end{equation}
By adjoining $F^A \comp G^B$ to the left and then $A \times B$ to the right we obtain maps
\begin{equation}
\label{dpair}
( F^A \comp G^B ) \otimes (A \times B) \xto{\rho^{F,G}_{A,B}} F \comp G 
\quad \text{and} \quad
F^A \comp G^B \xto{\Theta_{A,B}^{F,G}} (F \comp G)^{A\times B}.
\end{equation}
They are natural in $A,B,F$ and $G$ and we will frequently omit these ``decorations'' from $\la^{F,G}_{A,B}$.
By adjunction, \eqref{defrhoadjoint} is equal to
\begin{equation}\label{adj-theta-twice}
A \times B \xto{\ev_{A \times B}^\#} \bigl\{(F\comp G)^{A \times B}, F \comp G\bigr\}
\xto{(\Theta_{A,B}^{F,G})^*} \bigl\{F^A \comp G^B, F \comp G\bigr\}.
\end{equation}
We claim that the natural transformations $\Theta$ are associative in the sense that
\begin{equation}\label{theta-associative}
\begin{matrix}
\xymatrix@1{
[A,F] \comp [B,G] \comp [C,H] 
\ar[rr]^{\Theta^{F,G}_{A,B} \comp [C,H]}
\ar[d]_{[A,F] \comp \Theta^{G,H}_{B,C}}
& &
[A \times B, F \comp G] \comp [C,H] 
\ar[d]^{\Theta_{A \times B,C}^{F \comp G,H}} 
\\
[A,F] \comp [B \times C,G \comp H] 
\ar[rr]_{\Theta_{A, B \times C}^{F , G \comp H}}
&&
[A \times B \times C , F \comp G \comp H]
}
\end{matrix}
\end{equation}
commutes. 
Indeed, consider first the composition of the arrows at the top and right.
After adjoining $A \times B \times C$ to the left and $F^A \comp G^B \comp H^C$ to the right, 
the adjunctions of \eqref{defrhoadjoint} and \eqref{dpair} show that
this map becomes the composition of the arrows on the left and bottom of the following diagram.
$$
\xymatrix{
A\times B \times C
\ar[r]^(0.4){\ev_A^\# \times \ev_B^\# \times \ev_C^\#}  
\ar[d]^{\ev_{A\times B}^\# \times \ev_C^\#} 
\ar@/_6.5pc/[dd]_(0.7){(\rho^{F \diamond G,H}_{A \times B,C})^\#}
&
\{ F^A,F\} \times \{ G^B,G\} \times \{ H^C,H\}
\ar[d]_{\Diamond \times \Id}
\ar@/^6.5pc/[dd]^(0.7){\Diamond}
\\
\{ (F\comp G)^{A \times B},F\comp G\} \times \{H^C,H\}
\ar[r]^{\Theta^*\times \Id}
\ar[d]^{\Diamond} 
&
\{ F^A \comp G^B,F\comp G\} \times \{H^C,H\}
\ar[d]_{\Diamond}
\\
\{(F\comp G)^{A\times B} \comp H^C, F\comp G\comp H\}
\ar[r]_{(\Theta\comp H^C)^*}
&
\{ F^A\comp G^B\comp H^C,F\comp G\comp H\}
}
$$
The first square of this diagram commutes by \eqref{adj-theta-twice} and the second since $\Diamond$ is functorial.
Using a similar argument one shows that the adjoint of the composition of the arrows at the left and bottom of  \eqref{theta-associative} is also equal to the composition of the arrows at the top and right of the diagram above.
%
It follows that \eqref{theta-associative} commutes.

Now consider a small CW-category $C$ (\ref{toplogicalcategories}) and a functor $F \colon C \to \E$.
For any $A \colon C \to \Top$ we define
\[
[A,F]_C = \int_{c \in C} [A(c),F(c)] \in \E \quad \text{(cf. Definition~\ref{holimdef}).}
\]
The assignment $F \mapsto [A,F]_C$ is a continuous functor because the inverse limit functor is.

Given functors $F_i \colon C_i \to \E$ ($i=1,2$), there is a functor
\[
F_1 \extcomp F_2 \colon C_1 \times C_2  \to \E, \qquad 
(c_1,c_2) \mapsto F_1(c_1) \comp F_2(c_2).
\]
The transformation $\extcomp$ is natural and associative because the monoidal operation $\comp$ is associative. This construction generalizes \eqref{def-exttimes}.

By taking ends, the natural transformations $\Theta$ of \eqref{dpair} now give rise to natural transformations 
\begin{equation}\label{def-exttheta}
\Theta \colon [A_1,F_1]_{C_1} \comp [A_2,F_2]_{C_2} \to [A_1 \exttimes A_2 , F_1 \extcomp F_2]_{C_1 \times C_2}.
\end{equation}
where $A_i \in C_i^*$ and $F_i \colon C_i \to \E$ for $i=1,2$. The associativity \eqref{theta-associative} and the naturality of limits imply a similar associativity:
\begin{equation}\label{exttheta-associative}
\begin{matrix}
\xymatrix{
[A_1,F_1]_{C_1} \comp [A_2,F_2]_{C_2} \comp [A_3,F_3]_{C_3}
\ar[r]^-{\Theta \comp \Id }
\ar[d]_{\Id \comp \Theta}
&
[A_1 \exttimes A_2, F_1 \extcomp F_2]_{C_1 \times C_2} \comp [A_3,F_3]_{C_3}
\ar[d]^{\Theta}
\\
[A_1,F_1]_{C_1} \comp [A_2 \exttimes A_3, F_2 \extcomp F_3]_{C_2 \times C_3} \ar[r]_-{\Theta} 
&
[A_1 \exttimes A_2 \exttimes A_3, F_1 \extcomp F_2 \extcomp F_3]_{C_1 \times C_2 \times C_3}
}
\end{matrix}
\end{equation}

We will now fix an $A_\infty$-operad $A$ with category of operators $\Dplus=\Deplus(A)$ (Definition~\ref{defdeltaA}).
The inclusion of the full subcategory of the non-empty sets is denoted $J \colon D \to \Dplus$.
We recall from Definition~\ref{def-k-completion} that an $A$-algebra $K$ gives rise to a monoidal functor 
$R=R_K \colon \Dplus \to \E$ such that $R_K(*)=K$.
By Definition \ref{def-completion} and \eqref{TEDasLan}, $\hat{K}= [ED,J^*R_K]_D = [\TED,R_K]_{\Dplus}$.

Since $R$ is monoidal, for every finite ordered set $S$ (see Definition~\ref{ord-sum-functors}),
\begin{equation}\label{u-star-r}
\underbrace{R \extcomp \cdots \extcomp R}_{|S| \text{ times}} = u(S)^*(R).
\end{equation}

\begin{defn}\label{def-RS}
For any finite ordered set $S$, define
$$
R^{(S)} = \Bigl[\bigexttimes_{s \in S} \TED, \bigextcomp_{s \in S} R\Bigr]_{\Dplus^S} \in \E
$$
\end{defn}

In particular, $R^{(1)}=\hat{K}$.
The natural transformations $\Theta$ of \eqref{def-exttheta} give rise to natural transformations
$$
R^{(S)} \comp R^{(T)} \to R^{(S \sqcup T)}.
$$
It follows from \eqref{exttheta-associative} that for every monotonic $\vp \colon S \to T$ there results a well-defined natural transformation
\begin{equation}\label{def-beta}
\be_\vp \colon \bigcomp_{t \in T} R^{(\vp^{-1}t)} \to R^{(S)}, \qquad \text{(see Notation \ref{operadnotation}).}
\end{equation}
It also follows from \eqref{exttheta-associative} that the transformations $\be$ are associative in the sense
that for any $S \xto{\vp} T \xto{\psi} W$ the following square commutes:
\begin{equation}\label{be-associative}
\xymatrix@C=3cm{
{\underset{w \in W}{\bigcomp}} \ \underset{t \in \psi^{-1}w}{\bigcomp} R^{(\vp^{-1}t)}
\ar[r]^{\underset{w \in W}{\bigcomp} \be_{\vp|_{\vp^{-1}\psi^{-1}w}}}
\ar@{=}[d]
&
\underset{w \in W}{\bigcomp} R^{(\vp^{-1}\psi^{-1}w)}
\ar[d]^{\be_{\psi\circ\vp}}
\\
{\underset{{t \in T}}{\bigcomp}} R^{(\vp^{-1}t)}
\ar[r]_{\be_{\vp}}
&
R^{(S)}
}
\end{equation}
We see from \eqref{u-star-r}, from the adjunction of $u(S)^*$ and $\Lan_{u(S)}$, and from Definition~\ref{def-epsilonS} that
\begin{equation}\label{RS-simple}
R^{(S)} = [\epsilon(S),R]_{\Dplus}.
\end{equation}
Thus, there is a continuous map
$$
\map_{\Dplus^*}(\epsilon(S),\epsilon(T)) \xto{ \ f \mapsto f^* \ } \map(R^{(S)}, R^{(T)}).
$$

\begin{lem}\label{belan}
Let $\vp \colon S \to T$ be a monotonic map and consider $f_t \colon \TED \to \epsilon(\vp^{-1}t)$  for all 
$t \in T$, cf. Notation \ref{operadnotation}. 
Then the following square commutes:
\[
\xymatrix@C=2cm{
{\bigcomp_{t \in T}} R^{(\vp^{-1}t)} 
\ar[r]^{\bigcomp_{t \in T} (f_t^*)}
\ar[d]_{\be_\vp}
&
\bigcomp_{t \in T} R^{(1)}
\ar[d]^{\be_{\id_T}}
\\
R^{(S)}
\ar[r]^{(\Lan_{u(T)} (\bigexttimes_t f_t))^*} 
&
R^{(T)}
}
\]
\end{lem}

\begin{proof}
Immediate from the naturality of $\Theta$ together with 
\eqref{RS-simple} and \eqref{epsilon-of-sqcup}.
\end{proof}

\begin{proof}[Proof of Theorem~\ref{extendedoperadaction}]
For every finite ordered set $S$ define
\[
\si_S \colon \AHP(S) \to \Map(\bigcomp_S R^{(1)}, R^{(1)})
\]
by the composition
\begin{equation}\label{def-ahp-action}
\map_{\Dplus^*}(\TED,\epsilon(S)) \xto{f \mapsto f^*} \Map(R^{(S)},R^{(1)}) \xto{\be_{\id_S}^*}
\Map(\bigcomp_S R^{(1)},R^{(1)}).
\end{equation}
We claim that these maps form a morphism of operads $\AHP \to \opEND(\hat{K})$.
First, observe that when $S=*$ then $\be_{\id_S} = \id_{R^{(1)}}$ and therefore 
$\si(\id_{\TED})=\id_{R^{(1)}}$.
Thus, $\si$ respects the identity elements of the operads.

It remains to prove that $\si$ respects the composition operation in both operads, which we denote by $\bullet$.
Consider some $\vp \colon S \to T$ and elements
\begin{eqnarray*}
&& f \in \map(\TED,\epsilon(T)) = \AHP(T), \qquad \text{ and} \\
&& g_t \in \map(\TED,\epsilon(\vp^{-1}t))=\AHP(\vp^{-1}t) \qquad (t \in T).
\end{eqnarray*}
By construction (Definition~\ref{construction-AHP}),  $f \bullet (g_t) = (\Lan_{u(T)} \bigexttimes_t g_t) \circ f$.
Now, by the definition of $\si$, the monoidality of $\comp$, Lemma~\ref{belan}, and \eqref{be-associative},
\begin{align*}
\si(f) \bullet \bigl(\si(g_t)\bigr)_t  = &
(f^* \circ \be_{\id_T}) \circ \Bigl(\bigcomp_{t \in T} (g_t^* \circ \be_{\id_{\vp^{-1}t}})\Bigr) \\ 
= & (f^* \circ \be_{\Id_T}) \circ \Bigl(\bigcomp_{t \in T} g_t^*\Bigr)  \circ \Bigl(\bigcomp_{t \in T} \be_{\id_{\vp^{-1}t}}\Bigr)  
\\
= & f^* \circ \Bigl(\Lan_{u(T)}\bigexttimes_t g_t\Bigr)^* \circ \be_\vp \circ \Bigl(\bigcomp_{t \in T} \be_{\id_{\vp^{-1}t}}\Bigr)\\
= &
\bigl(f \bullet (g_t)_t\bigr)^* \circ \be_{\id_S} = \si\bigl(f \bullet (g_t)_t\bigr). 
\end{align*}
That is, $\si$ respects the composition operations of the operads.
\end{proof}

%
%
%
%
\section{The restricted operad action}
\label{restriced-operad-section}

In Section~\ref{action-on-completion} we showed how $\AHP$ acts on the completion $\hat{K}$ of any $A$-algebra $K$.
The goal of this section is to construct an operad $\hat{A}$ which is $A_\infty$ if $A$ is, together with morphisms of
operads $\AHP \leftarrow \hat{A} \to A$ which make the natural map $\hat{K} \to K$ defined in \eqref{khat-to-k} into
a morphism of $\hat{A}$-algebras.

Recall that for any $\Top$-category $C$ we write $y_c$ for the corepresentable functor $\map_C(c,-)$.
An element $g \in \map_C(c,c')$ gives rise to a natural transformation $g^* \colon y_{c'} \to y_c$ which we also
denote by $y_g$.

Fix an operad $A$ and set $\Dplus = \Deplus(A)$, see Definition~\ref{defdeltaA}.
Following Notation~\ref{operadnotation}, given a monotonic map $\vp \colon S \to T$, we obtain functors $y_{\vp^{-1}t}$.
By (\ref{epsilon-of-sqcup}, \ref{lan-exttimes-corep}), $\Lan_{u(T)} (\bigexttimes_t y_{\vp^{-1}t}) = y_{S}$.
More generally, for $\psi \colon S' \to S$ we have 
\begin{equation}\label{lkan-yoneda}
\Lan_{u(T)}\Bigl(\bigexttimes_t g_t^*\Bigr)= \Bigl(\bigsqcup_t g_t\Bigr)^* \colon y_{S} \to y_{S'}
   \qquad \text{where } g_t \in \map_{\Dplus}(\psi^{-1}\vp^{-1}t,\vp^{-1}t).
\end{equation}
Since $\E_0 \Dplus = \coprod_{T \in \Dplus} y_T$, we obtain from $\E_0\Dplus \subseteq E\Dplus$ a natural transformation
\[
\kappa_\pnt \colon y_\pnt \to \TED.
\]
More generally, from \eqref{lkan-yoneda} and Definition~\ref{def-epsilonS} we obtain a natural transformation
\begin{equation}\label{defkappaS}
\kappa_S \colon y_S \to \epsilon(S) \qquad \text{defined by } \kappa_S = \Lan_{u(S)} \Bigl(\bigexttimes_S \kappa_\pnt\Bigr) \qquad\text{($S \in \Dplus$)}.
\end{equation}
For every finite ordered set $S$, we define maps
\begin{eqnarray*}
&& \rho \colon \AHP(S) = \map_{\Dplus^*}(\TED,\epsilon(S)) \xto{(\kappa_\pnt)^*} \map_{\Dplus}\bigl(y_\pnt,\epsilon(S)\bigr), \\
&& \zeta \colon A(S) = \map_{\Dplus}(S,*) = \map_{\Dplus^*}(y_\pnt,y_S) \xto{(\kappa_S)_*} \map_{\Dplus^*}\bigl(y_\pnt,\epsilon(S)\bigr).
\end{eqnarray*}

\begin{defn}\label{def-a-hat}
For any finite ordered set $S$, define the space $\hat{A}(S)$ as the pullback 
$$
\xymatrix@1{
{\hat{A}}(S)
\ar[r]^{\bar{\zeta}}
\ar[d]_{\bar{\rho}} 
&
\AHP(S)
\ar[d]^{\rho}
\\
A(S) \ar[r]_-{\zeta} &
\map_{\Dplus^*}\bigl(y_\pnt,\epsilon(S)\bigr)
}
$$
\end{defn}

Thus, $\hat{A}(S)$ is the subspace of $\AHP(S) \times A(S)$ consisting of the pairs $(f,g)$ with $f \colon \TED \to \epsilon(S)$ and $g \in A(S)=\map_{\Dplus}(S,*)$ satisfying
\begin{equation}\label{who-is-in-hatA}
f \circ \kappa_\pnt = \kappa_S \circ y_g.
\end{equation}

We equip $\AHP \times A$ with the product operad structure.

\begin{thm}
\label{thm-hatA-contractible}
The spaces $\hat{A}(S)$ form a sub-operad of $\AHP \times A$.
In particular the projections $\hat{A} \to \AHP$ and $\hat{A} \to A$ are morphisms of operads.
Moreover, if $A$ is an $A_\infty$-operad then 
$\hat{A}(S)$ is weakly contractible for all $S$ and if $A(\emptyset)=*$ then $\hat{A}(\emptyset)=*$.
\end{thm}

\begin{proof}
We first note that if $A(\emptyset)=*$ then $\epsilon(\emptyset)=*$ 
by Lemma \ref{ordsumeasy}\eqref{uemptyset} and Definition~\ref{def-epsilonS}.
Therefore all three spaces defining the pullback square of Definition~\ref{def-a-hat} are points, whence
$\hat{A}(\emptyset)=*$.

Now assume that $A$ is an $A_\infty$-operad.
Consider the inclusion $y_\pnt \subseteq \E_0 \Dplus = \coprod_{T \in \Dplus} y_T$.
By Lemma~\ref{ECcofibrant} there is a composite fibration
$$
\map_{\Dplus^*}(\TED,\epsilon(S)) \to \map_{\Dplus^*}(\E_0 \Dplus,\epsilon(S)) \to \map_{\Dplus}(y_\pnt,\epsilon(S))
=\epsilon(S)(*)
$$
because the second map is isomorphic to the fibration 
$\prod_{T \in \Dplus} \epsilon(S)(T) \to \epsilon(S)(*)$ by Yoneda's Lemma.
In the proof of Theorem~\ref{ahpcontractible} we have shown that $\epsilon(S)(*)$ is weakly contractible and deduced
that the spaces $\AHP(S)$ are weakly contractible.
It now follows that the map $\rho$ in the pullback square in Definition~\ref{def-a-hat} is an equivalence as well as a fibration.
Since a pullback of a trivial fibration of spaces is a trivial fibration, we deduce that $\hat{A}(S) \to A(S)\simeq *$ is
an equivalence.

It remains to show that the spaces $\hat{A}(S)$ form a sub-operad of $\AHP \times A$.
We will let $\bullet$ denote the operad composition in both operads $\AHP$ and $A$.
Fix a monotonic map $\vp \colon S \to T$ and consider $(f,g) \in \hat{A}(T)$ and $(f_t,g_t) \in \hat{A}(\vp^{-1}t)$.
Set $f'=f \bullet(f_t)$ and $g'=g \bullet (g_t)$.
We now claim that $(f',g') \in \hat{A}(S)$. First, we note that by \eqref{Lkan-product-functors}, Lemma~\ref{ordsumeasy}\eqref{uS-in-steps}, the definition of $\kappa_{\vp^{-1}t}$, and \cite[Theorem~4.47]{kelly:enriched},
\[
\Lan_{u(T)} \Bigl(\bigexttimes_{t \in T} \kappa_{\vp^{-1}t}\Bigr) = \Lan_{u(T)} \biggl(\bigexttimes_{t \in T} \Lan_{u(\vp^{-1}t)}\Bigl(\bigexttimes_{\vp^{-1}t} \kappa_\pnt\Bigr)\biggr) = \Lan_{u(S)}\Bigl(\bigexttimes_S \kappa_\pnt\Bigr) = \kappa_S.
\]
We also observe that by Definition~\ref{defdeltaA} of the category $\Dplus$ and \eqref{lkan-yoneda},
\begin{equation}\label{dplus-yoneda-compose}
\Lan_{u(T)}\Bigl(\bigexttimes_{t \in T} g_t^*\Bigr) \circ g^* = \Bigl(\bigsqcup_t g_t\Bigr)^* \circ g^* = \bigl(g \bullet (g_t)\bigr)^*.
\end{equation}
It follows that $(f',g') \in \hat{A}(S)$ because \eqref{who-is-in-hatA} applies to $(f,g)$ and $(f_t,g_t)$, hence
\begin{align*}
f' \circ \kappa_\pnt =& 
\Lan_{u(T)}\Bigl(\bigexttimes_T f_t\Bigr) \circ f \circ \kappa_\pnt = 
\Lan_{u(T)}\Bigl(\bigexttimes_T f_t\Bigr) \circ \kappa_T \circ g^* \\
\underset{\eqref{defkappaS}}{=} &
\Lan_{u(T)}\Bigl(\bigexttimes_T (f_t \circ \kappa_\pnt)\Bigr) \circ g^* = 
\Lan_{u(T)}\Bigl(\bigexttimes_T (\kappa_{\vp^{-1}t} \circ g_t^*)\Bigr) \circ g^*
\\
= &
\Lan_{u(T)}\Bigl(\bigexttimes_T \kappa_{\vp^{-1}t}\Bigr) \circ \Lan_{u(T)} \Bigl(\bigexttimes_T g_t^*\Bigr) \circ g^* \underset{\eqref{dplus-yoneda-compose}}{=} 
\kappa_S \circ (g \bullet (g_t))^* = \kappa_S \circ (g')^*.
\end{align*}
Thus, $\hat{A}$ is stable under the composition law of $\AHP \times A$.
It remains to show that $\hat{A}$ contains the identity element $(\id_{\TED},\iota)$ of $\AHP \times A$.
Indeed, $\iota$ is the identity morphism in $\map_{\Dplus}(*,*)$ and therefore $\id_{\TED} \circ \kappa_\pnt = \kappa_\pnt \circ \iota^*$.
Thus, $(\id_{\TED},\iota) \in \hat{A}(*)$ by \eqref{who-is-in-hatA}.
\end{proof}

If $K$ is an $A$-algebra, $\kappa_\pnt$ gives rise to a natural transformation
\begin{equation}\label{khat-to-k}
\hat{K} = [\TED,R_K]_{\Dplus} \xto{\quad \tau=(\kappa_\pnt)^* \quad } [y_*,R_K]_{\Dplus} = R_K(*) = K.
\end{equation}

\begin{proof}[Proof of Theorem~\ref{completion-theorem}]
We consider the morphisms of operads $\hat{A} \to \AHP$ and $\hat{A} \to A$ obtained from Theorem~\ref{thm-hatA-contractible}.
Together with Theorem~\ref{extendedoperadaction} we obtain a morphism of operads
$$
\hat{A} \to \AHP \to \opEND(\hat{K}),
$$
that is, $\hat{K}$ is an $\hat{A}$-algebra.
To prove that $\tau \colon \hat{K} \to K$ is a morphism of $\hat{A}$-algebras we need to prove that that the following diagram is commutative:
\begin{equation}\label{monad-map:1}
\begin{matrix}
\xymatrix{
{\hat{A}}(S) \ar[r]^-{\si_{\hat{K}}} \ar[d]_{\bar{\rho}} &
\Map\Bigl(\bigcomp\limits_S R^{(1)},R^{(1)}\Bigr) \ar[r]^{\tau_*} &
\Map\Bigl(\bigcomp R^{(1)},K\Bigr) 
\\
A(S) \ar[rr]_{\si_K} & &
\Map\Bigl(\bigcomp K, K\Bigr) \ar[u]_{(\bigcomp \tau)^*} 
}
\end{matrix}
\end{equation}
By Definition~\ref{def-k-completion}, $K$ gives rise to a monoidal functor 
$R = R_K \colon \Deplus(A) \to \E$ such that $R_K(*)=K$ and $R^{(1)}=\hat{K}$ \eqref{def-RS}.
The naturality of the transformation $\Theta$ \eqref{def-exttheta} and the definition of $\be$ \eqref{def-beta} imply that the following diagram commutes:
\begin{equation}\label{trapezoiddiagram}
\begin{matrix}
\xymatrix{
{\bigcomp_S [\TED,R]_{\Dplus}} \ar[r]^{\bigcomp_S \tau} \ar[d]_{\be_{\id_S}} &
\bigcomp_S [y_*,R]_{\Dplus}
\ar[d]^{\Theta}
\ar[r]^{\si_K(g)} &
[y_\pnt,R]_{\Dplus}
\\
[\epsilon(S),R]_{\Dplus} \ar[r]_{(\kappa_S)^*} &
[y_S,R]_{\Dplus}
\ar[ur]_{(y_g)^*}
}
\end{matrix}
\end{equation}
By Yoneda and because $(y_g)^*$ is isomorphic to $\si_K(g)=R(g) \colon R(S) \to R(*)$,
\[
\bigcomp_S \bigl[y_\pnt,R\bigr]_{\Dplus} = \bigcomp_S R(*) = R(S) = [y_S,R]_{\Dplus}.
\]

The image of $(f,g) \in \hat{A}(S) \subseteq \AHP(S) \times A(S)$ under the top map of Diagram~\eqref{monad-map:1} is the natural transformation 
$$
\bigcomp_S R^{(1)} \xto{\be_{\id_S}} R^{(S)} \xto{f^*} R^{(1)} \xto{\tau=(\kappa_\pnt)^*} [y_\pnt,R]=K.
$$
Since $f \circ \kappa_\pnt = \kappa_S \circ g^*$ by \eqref{who-is-in-hatA}, the commutativity of \eqref{trapezoiddiagram} implies
\[
\tau \circ f^* \circ \be_{\id_S} = (\kappa_S \circ g^*)^* \circ \be_{\id_S} = (y_g)^* \circ \kappa_S^* \circ \be_{\id_S} = 
\Bigl(\bigcomp_S \tau\Bigr) \circ \si_K(g).
\]
This shows that \eqref{monad-map:1} commutes, hence $\hat{K} \to K$ is a morphism of $\hat{A}$-algebras.

To resolve the difficulty that $\hat{A}$ need not consist of CW-complexes, we apply the map
$|\Sing(\hat{A})| \to \hat{A}$. This is a morphism of operads by naturality, since both $\Sing$ and $| - |$
commute with products, and $|\Sing(*)|=*$.
\end{proof}


\end{document}